\documentclass[12pt,leqno]{article}
\usepackage{amsmath,amssymb}
\newcommand{\be}{\begin{equation}}
\newcommand{\ee}{\end{equation}}
\newcommand{\bea}{\begin{eqnarray}}
\newcommand{\eea}{\end{eqnarray}}
\newcommand{\bean}{\begin{eqnarray*}}
\newcommand{\eean}{\end{eqnarray*}}

\newcommand{\Q}{\mathbb{Q}}
\newcommand{\R}{\mathbb{R}}
\newcommand{\C}{\mathbb{C}}
\newcommand{\zz}{\mathbb{Z}}
\newcommand{\go}{\omega}
\newcommand{\ga}{\alpha}
\newcommand{\gb}{\beta}

\newenvironment{caja}[1]{\begin{trivlist}
       \item[\hskip \labelsep {\bfseries #1}]}{\end{trivlist}}

\begin{document}

\title{On the Hodge Conjecture for products of certain surfaces}
\author{Jos\'e J. Ram\'on Mar\'i \\Humboldt Uni. zu Berlin\\ Institut f\"ur Mathematik \\ Unter den Linden, 6  \\ D - 10099 Berlin \\Germany\\E-mail: akalmahakal@hotmail.com, \\jjramon@mathematik.hu-berlin.de}

\maketitle
\newtheorem{theorem}{Theorem}[section]
\newtheorem{prop}[theorem]{Proposition}
\newtheorem{corollary}[theorem]{Corollary}
\newtheorem{lemma}[theorem]{Lemma}
\newtheorem{obs}[theorem]{Remark}
\abstract{In this paper we prove the Hodge conjecture for arbitrary products of surfaces, $S_1 \times \cdots \times S_n$ such that $q(S_i)=2, p_g(S_i)=1$.  We also prove the Hodge conjecture for arbitrary self-products of a K3 surface $X$ such that the field $E=End_{hg} T(X)$ is CM.}

\vspace*{1cm}

\noindent
\section*{0  Notation and preliminaries}

  Unless otherwise stated, we use the terms \textit{curve} and \textit{surface} to denote smooth projective curves and surfaces, resp.  The term $p_g(S)=h^{2,0}(S)$ is called the \textit{geometric genus} of $S$, and $q(S)=h^{1,0}(S)=\mbox{dim }Alb(S)$ is known as the \textit{irregularity} of $S$.  For any complex projective manifold $X$, $H^k(X)$ will denote the group $H^k(X,\Q)$ regarded as a (rational) Hodge structure of (pure) weight $k$.  All Hodge structures appearing in this paper are rational and pure ~\cite{DeligneK3}; as usual, a \textit{Hodge cycle} (of \textit{codimension} $p$) or Hodge class of a Hodge structure $V$ is an element $v \in V_\C^{p,p} \cap V$.  We denote the subspace of Hodge cycles of $V$ by ${\cal H}(V)$, and also ${\cal H}^p(X)={\cal H}(H^{2p}(X))$ for $X$ a smooth projective variety; consequently, ${\cal H}(X)=\oplus_{i=0}^{\mathrm{dim}(X)} {\cal H}^i(X)$ will denote the \textit{Hodge ring}, or ring of Hodge classes of $X$.

    We define the (rational) \textit{transcendental lattice} $T(S)$ of a surface $S$ by the following orthogonal decomposition
\be \label{trasc}
H^2(S)=T(S) \oplus NS(S)_\Q
\ee
with respect to the cup-product.  The cup-product induces, after a change of sign, a polarisation of the Hodge structure $T(S)$ ~\cite{DeligneK3}.

  For $V$ and $W$ two (pure) Hodge structures of the same weight, we denote $Hom_{hg}(V,W)$ to be the space of linear maps from $V$ to $W$ respecting the Hodge structures.  For an introduction see ~\cite{DeligneK3}, ~\cite{vG}.

  For a Hodge structure $V$ as above we define the \textbf{Hodge group} of $V$, $Hg(V)$ to be the minimal $\Q$-defined algebraic subgroup of $GL(V)$ such that $h(U(1)) \subset Hg(V)_\R$; here $h$ is the representation corresponding to the Hodge bigraduation as in ~\cite{DeligneK3}.  The following is basic in this paper:

\begin{prop}\label{mtss} ~\cite{DeligneK3} ~\cite{DMOS}  Let $V$ be a polarisable Hodge structure.  Then $Hg(V)$ is reductive.  As a result, the category of polarisable Hodge structures is semisimple abelian.
\end{prop}

  For an comprehensive survey on the Hodge conjecture for abelian varieties, as well as a detailed introduction on the Hodge group $Hg(A)$, we refer the reader to ~\cite{Lewis} Appendix B.
\begin{caja}{Acknowledgments:}  I thank Prof. Bert van Geemen for a very useful and generous tutorial on the subject.  I wish to express my gratitude to Prof. B.J. Totaro for valuable suggestions on the writing of this paper.  I am grateful to the EU Research Training Network 'Arithmetic Algebraic Geometry' for their financial support during my PhD.  The warm hospitality of the Isaac Newton Institute is gratefully acknowledged.  I finally wish to thank Prof. G. E. Welters for introducing me to this subject.
\end{caja}
\section{Introduction}
  The purpose of this article is to prove the Hodge conjecture in two different situations of product of surfaces.  The first one is the product $S_1 \times \cdots \times S_n$, where $q(S_i)=2$ and $p_g(S_i)=1$. (It turns out that these surfaces are birationally equivalent to abelian or elliptic isotrivial surfaces).  This result  generalises the Main Theorem in the author's PhD thesis ~\cite{T}.

  The other case we consider is the following: take a K3 surface $X$; then the transcendental lattice $T(X)$ is irreducible, and its endomorphism algebra is a number field $E=End_{hg}T(X)$, which can be either CM or totally real (~\cite{Zarhin} 1.5).  We will prove that the Hodge conjecture for arbitrary powers of $X$ follows from the Hodge conjecture for $X \times X$.  In the case when $E$ is CM, we use results of Mukai ~\cite{Mukai}, together with an elementary lemma, to prove the Hodge conjecture for $X \times X$, and establish the result for $X^n$ for all $n$ by using invariant theory (see for instance ~\cite{Ribet} for similar arguments).
\section{Surfaces $S$ with $p_g=1, q=2$}
  This section is devoted to understanding the geometry of surfaces with $p_g=1, q=2$.

\begin{prop}\label{uno} Let $S$ be a minimal surface with $p_g=1, q=2$.  If $S$ is not abelian, then $S$ is of the form 
$$S= (C' \times E')/G$$ where $C'$ is a curve, $E'$ is an elliptic curve and $G$ acts faithfully on both components.  
\end{prop}
{\bf Proof:} One has $\chi({\cal O}_S)=0=1-q+p_g$.  From Enriques' classification we see that $S$ is non-ruled, and $K^2 \geq 0$.  Also $e(S)\geq 0$ (see ~\cite{Beauville} Th. X.4), and by Noether's formula we get $0=e(S)$, \textit{i.e.} $b_2(S)=6$, and so therefore $K^2=0$, which yields $S$ elliptic. Finally, by ~\cite{Serrano} or ~\cite{Beauville} Exs. VI.22(4), VIII.22., we see that $S=(C' \times E')/G$ is a finite \'etale quotient such that $g(E')=1$, and the proof is thus complete. $\blacksquare$


  All the statements concerning motives are, unless otherwise stated, considered in the category of Chow motives modulo homological equivalence.  We refer the reader to  ~\cite{Scholl} for the basic notations and language.
\begin{prop}{(Murre)}~\cite{Scholl} Let $X$ be a surface.  Then there exists a decomposition $h(X)= \oplus_{i=0}^{4} h^i(X)$; \textit{i.e.} a Chow-K\"unneth decomposition exists in the case of surfaces (in fact, modulo rational equivalence).
\end{prop}
\begin{obs} From the above and the standard conjectures for abelian varieties ~\cite{Kleiman}, it follows that the Hodge classes on $X\times X$ inducing the projectors $H^\bullet(X)\twoheadrightarrow H^i(X) \subset H^\bullet(X)$ on a variety $X$ which is a product of surfaces and abelian varieties are all algebraic, and thus $X$ admits a decomposition $h(X) \simeq \oplus_{i=0}^{2 dim(X)} h^i(X)$ modulo homological equivalence.   This result will be used throughout.
\end{obs}

\begin{prop}\label{motivos}
  Let $S$ be a minimal, not abelian surface such that $p_g=1, q=2$.  Notations being as in Proposition \ref{uno}, the following cases hold:
\begin{enumerate}
\item[(a)] either $g(E'/G)=1$ and the Albanese map $a$ induces an isomorphism $h(S) \simeq h(\mbox{Alb }S)$, or 
\item[(b)] $g(E'/G)=1$ and the Albanese map $a$ sends $S$ onto a curve $B$.  It turns out that $B=C'/G$ and the Albanese fibration $$S=(C'\times E')/G \rightarrow a(S)=B=C'/G$$ is the canonical projection.
\end{enumerate}
\end{prop}

\begin{caja}{Proof:}
The following argument holds in both cases ~\cite{Freitag} ~\cite{Serrano}: $H^1((C' \times E')/G)=H^1(C'/G)\oplus H^1(E'/G)$.
Since $q(S)=\frac{1}{2}b_1(S)$, we have 
\be \label{qs}
q(S)=g(C'/G)+g(E'/G)
\ee and so the following cases are possible.
\begin{enumerate}
\item[(a)]$g(E'/G)=1$, $g(C'/G)=1$.  In this case $G$ acts on $E'$ by translations, and $A=C'/G \times  E'/G$ is an abelian surface; the natural map 
$$\phi:S=(C' \times E')/G \rightarrow C'/G \times E'/G=A$$
yields an isomorphism on $H^1$ by the above (and so on $H^3$); therefore $\mbox{Alb }S \sim A$, whence $h^1(S) \simeq h^1(A)$.  On $H^2$, the following holds:
$$H^2(S)=H^2(C'\times E')^G= H^2(C'/G \times E'/G),$$
for $G$ acts freely on $C' \times E'$ and trivially on $H^\bullet(E')$; this proves that $\phi|H^2$ is an isomorphism, thus establishing the result.
\item[(b)] Let $B=C'/G$.  In this case we have $g(E'/G)=0$ and $g(B)=2$.  The natural map
$$p:S=(C' \times E')/G \rightarrow C'/G=B$$
 satisfies $q(S)=g(B)$ by Formula (\ref{qs}), and therefore coincides with the Albanese fibration ~\cite{Beauville}; see also ~\cite{Ueno} Ch. 9.  $\blacksquare$
\end{enumerate}
\end{caja}

\subsection{The case $g(E'/G)=0$}

  Let $S$ satisfy case \textit{(b)} of Proposition \ref{motivos}, and let $H\subset G$ be the subgroup of translations on $E'$.  Since $H-\{1\}$ coincides with the set of fixed-point-free transformations of $E'$ in $G$, we have a split exact sequence (we now fix a section $\sigma$)
$$1\rightarrow H\rightarrow G \leftrightarrows \zz_n\rightarrow 1,$$ 
where $\zz_n \hookrightarrow Aut_P(E')$ for $P$ fixed point of a generator $\phi$ of $\sigma(\zz_n)$.  Clearly $n \in \{2,3,4,6\}.$ 

  The following proposition is a reduction to the case $G=\zz_n, H=\{1\}$.

\begin{prop}\label{redu}
  Let $C=C'/H$, $E=E'/H$.  If the natural action $\mu$ of $\zz_n=G/H$ on $C \times E$ is \'etale, then the natural map
 $$\gb:S=(C'\times E')/G \rightarrow S'=(C \times E)/\zz_n$$ yields an isomorphism of motives $h(S) \simeq h(S')$.
\end{prop}
{\bf Proof:} The proof is similar to that of Proposition \ref{motivos}\textit{(b)}.

\subsection{$\mu$ is free}

 We suppose $g(E'/G)=0$, notations being as above.  We are going to prove that this case meets the hypotheses of Proposition \ref{redu}.
\begin{obs}\label{numeros}
 Consider the Hodge structure $$V=\left[H^1(C')\otimes H^1(E')\right]^G=\left[H^1(C)\otimes H^1(E)\right]^{\zz_n}.$$  Then

$$H^2(S)=V \oplus \Q(-1)^{\oplus 2}$$
and $V$ has Hodge numbers $\mbox{dim }V^{2,0}=1=\mbox{dim }V^{0,2},\mbox{dim }V^{1,1}=2.$
\end{obs}

Consider the action of $G/H$ on $JC$; let $\phi$ be a generator of $G/H$ such that $\phi^*|H^{1,0}(E)=\go$ where $\go=e^{2\pi i/n}$; let $Q_n(x)$ denote the cyclotomic polynomial of order $n$.  

\begin{theorem}\label{etale}
Let $P:=\mbox{ker }Q_n(\phi_*)^0 \subset JC$. Then $\mbox{dim }P=1$ for $n=2$ and $\mbox{dim }P=2$ for $n=3,4,6$.  The quotient map $C \rightarrow C/\zz_n=B$ is \'etale in all cases. 
\end{theorem}
{\bf Proof of Theorem  \ref{etale}:} Consider $V$ as above.  It is clear that
\be \label{igualdad}
V= \left(H^1(P)\otimes H^1(E)\right)^{\zz_n}.
\ee
In the case $n=2$, $\phi$ acts on both vector spaces as $-Id$, so $V=H^1(P)\otimes H^1(E)$, whence $\mbox{dim }P=1$ by inspection. 
For $n=3,4,6$, let $\chi$ be the character of $\zz_n$ such that $H^{1,0}(E)=\chi$; then $H^{1,0}(P)= a \chi \oplus b \overline{\chi}.$  Inspecting Hodge numbers as above and using Remark \ref{numeros} we find $a=b=1$, which in turn yields $\mbox{dim }P=2$.

  From the above we conclude that the action $\mu$ of $\zz_n$ on $C$ has no fixed points.  This follows from ~\cite{Schoen} Lemma 1.5; alternatively one can derive this result from several Riemann-Hurwitz type inequalities.
$\blacksquare$
\begin{corollary} The motive of a surface $S=(C' \times E')/G$ with $G \neq H$ is isomorphic to that of a surface $(C \times E)/ \zz_n$ with $H=\{1\}$.  In other words, the conclusion of Proposition \ref{redu} holds true always. 
\end{corollary}

\subsection{$h^2(S) \simeq h^2(A)$}

 We now consider $S$ as above, \textit{i.e. } with cyclic $G=\zz_n$, such that $B=C/G$ is a genus $2$ \'etale quotient, and find an abelian surface $A$ such that an isomorphism of Hodge structures $H^2(S) \cong H^2(A)$ holds.  The first step is to decompose $P$:
\begin{lemma}\label{spezza} The abelian surface $P$ above splits as $P \sim E_1 \times E_1$.
\end{lemma}
{\bf Proof:} Indeed, suppose that $P$ is simple.  Then $Hg(P \times E)=Hg(P)\times Hg(E)$  (due to F. Hazama; see e.g. ~\cite{Lewis}B.7.6.2; see also  ~\cite{MZ}), whence the Hodge structure $W=H^1(P) \otimes H^1(E)$ is irreducible (with $\mbox{dim }W^{2,0}=2$).  Hence $W$ cannot contain $V$, which contradicts our hypothesis.  Therefore  $P$ must split; using the $\zz_n$-decomposition of $H^1(P)$ from the Proof of Theorem \ref{etale} and an elementary argument we obtain $P \sim E_1 \times E_1$ for $E_1$ an elliptic curve, thereby completing the proof. $\blacksquare$

Let us get back to our $H^2(S)$.  We had by Formula (\ref{igualdad}) and Lemma \ref{spezza}
\be \label{inclfinal}
H^2(S)=\Q(-1)^2\oplus [H^1(P)\otimes H^1(E)]^{\zz_n} \subset \Q(-1)^2\oplus [H^1(E_1) \otimes H^1(E)]^2.
\ee

  Again, since the transcendental part of $H^1(E_1)\otimes H^1(E)$ has one-dimensional $(2,0)$-part (and is thus irreducible ~\cite{vG} ~\cite{Zarhin} ), by Formula (\ref{inclfinal})  $H^2(S)$ and $H^2(E_1 \times E)$ differ only by powers of the Tate Hodge structure, which implies $H^2(S)\cong H^2(E_1 \times E)$ by counting dimensions.  We have thus proven the following Proposition.

\begin{prop} Under the hypotheses of Proposition \ref{motivos}\textbf{(b)}, the abelian surface $A= E_1 \times E$ is such that  
$H^2(S) \cong H^2(A)$ (as Hodge structures).
\end{prop}$\blacksquare$

  We now proceed to construct an algebraic cycle inducing the described isomorphism.  The scheme is the following.  Choose a $\phi$-equivariant projection $u:JC \twoheadrightarrow P,$ and consider the correspondence $ \gb=(u_* \circ (alb_C)_*, id_E)\circ \pi^*$ from $S$ to $P \times E$, where $\pi:C \times E \rightarrow (C \times E)/\zz_n=S$ is the natural projection.  This correspondence from $S$ to $P \times E$ realises the inclusion in Formula (\ref{inclfinal}).  The final step in this construction will be to cook up a correspondence from $P \times E \sim E_1 \times E_1 \times E$ sending the image of $V$ onto $H^1(E_1) \otimes H^1(E)$ in $E_1 \times E$, which after composing can be easily extended to the sought-after isomorphism.

\begin{lemma}\label{corrfinal} Let $E_1, E_2$ be two elliptic curves.  For every Hodge substructure $V$ of $H^1(E_1 \times E_1)\otimes H^1(E_2)$ isomorphic to $H^1(E_1)\otimes H^1(E_2)$ there exists an algebraic correspondence $\ga$ from $E_1 \times E_1 \times E_2$ to $E_1 \times E_2$ such that $\ga_*V=H^1(E_1)\otimes H^1(E_2)$.
\end{lemma}
\begin{caja}{Proof of Lemma \ref{corrfinal}:} It suffices to prove that every Hodge correspondence between $H^1(E_1)\otimes H^1(E_2)$ and $H^1(E_1 \times E_1)\otimes H^1(E_2)$ is algebraic. This follows from the Hodge conjecture for products of elliptic curves, due to Imai ~\cite{Lewis} ~\cite{MZ} (see also Proposition \ref{abe} below.)
\end{caja}

  We are now ready to prove the following result:
\begin{theorem} With the assumptions of this Section, the motives $h^2(S)$ and $h^2(E_1 \times E)$ are isomorphic (modulo homological equivalence). 
\end{theorem}
{\bf Proof:}
  To conclude the proof, consider the correspondence $\gb$ above, which takes $V$ to its image inside $H^1(P)\otimes H^1(E)$ of $P \times E \sim E_1 \times E_1\times E$.  Choose a projection 
$$\ga:H^1(E_1\times E_1)\otimes H^1(E) \twoheadrightarrow H^1(E_1)\otimes H^1(E)$$  such that $\ga|V$ is a (Hodge) isomorphism.  $\ga$ is algebraic by Lemma \ref{corrfinal}, and so the composition $\ga \circ \gb$, also algebraic, yields the desired isomorphism.  
$\blacksquare$
\begin{obs} An explicit isomorphism could be obtained by fiddling with $\phi^*$ as an element of $M_2(End(E_1)\otimes End(E))$, without the use of Lemma \ref{corrfinal}.  We leave this to the reader.
\end{obs}
\subsection{The Hodge Conjecture for $S_1 \times \cdots \times S_m$, $p_g(S_i)=1, q(S_i)=2$}

We are going to prove the following theorem:
\begin{theorem}\label{principal} Let $S_i$ be surfaces such that $p_q(S_i)=1, q(S_i)=2$ ($S_i$ need not be minimal).  Then the Hodge conjecture holds for $S_1 \times \cdots \times S_m$.
\end{theorem}
\begin{obs}\label{ene}
 Let $S$ be a surface such that $p_g=1, q=2$.  In the former sections we have actually proven that the motive of such a surface (minimal or not) is generated (in the Tannakian sense, see ~\cite{DMOS}) by motives of abelian surfaces and elliptic curves.
\end{obs}
The following lemma follows easily from ~\cite{DMOS}(see also ~\cite{Lewis} Appendix B) and some linear algebra.
\begin{lemma} Let $A$ be an abelian variety of dimension $\geq 2$.  Then $Hg(H^2(A))=Hg(A)/\mu_2.$  In particular, for $A$ an abelian surface of simple CM type $(F, \Phi)$ one has $U_F(1)\simeq Hg(T(A))=Hg(A)/\mu_2=U_F(1)/\mu_2.$
\end{lemma}
$\blacksquare$

  There are two cases of a simple quartic CM field $F$ ~\cite{ST} according to the Galois group of its normal closure $N$ over $\Q$.   We now solve the case $F\neq N$.  In the case where $F=N$, there is essentially one CM-type, hence one isogeny type of abelian surfaces ~\cite{T}.  In the case when $F$ is Galois and $Gal(F|\Q)=V_4$ the corresponding abelian surface is non-simple ~\cite{ST}.  Therefore one is left with the cases $F=N$, $Gal(F|\Q)=C_4$ and $Gal(N|\Q)=D_{2,4}$.

\begin{lemma}\label{diedral} \begin{enumerate} \item Let $F$ be a simple CM field of degree $4$, Galois over $\Q$ (\textit{i.e.} $Gal(F|\Q)=C_4$).  Then $$End_{hg}(T(A))=F.$$

\item Let $F$ be a CM field of degree $4$ over $\Q$ such that its normal closure has Galois group $D_{2,4}$ and $A$ belong to the CM type $(F,\Phi),$ where $\Phi=\{ \sigma_1,\sigma_2 \}$. Then $End_{hg}T(A)=E,$ where $E \subset N$ is quartic CM and non-isomorphic to $F$.  As a result one has the following isomorphism of algebraic groups over $\Q$:
$$U_F(1)\simeq U_E(1),$$
with $E$ and $F$ non-isomorphic number fields.  In fact, if we write $F=F_0[\theta]$ where $\theta^2=-\ga$ for $\ga\in F_0$ totally positive, then $E\cap \R=\Q\left( \sqrt{N_{F_0|\Q}(\ga)}\right) \neq F_0$ (this follows from the condition on $Gal(N|\Q)$ and also from the uniqueness of $E$ up to automorphisms of $N|\Q$).  Thus $E$ and $F$ mutually determine each other.

 \end{enumerate}
\end{lemma}
\begin{caja}{Proof of Lemma \ref{diedral}:}

  One need only observe that the subfield $E$ of $\C$ spanned by the action of $F^\times$ on $T(A)$ (which can be read on $H^{2,0}(A)$) is quartic CM and not isomorphic to $F$.  Indeed, the homomorphism of abstract groups
$$\rho:F^\times \rightarrow GL(H^{2,0}(A))$$ is described by $x \mapsto \rho(x)=\sigma_1(x) \sigma_2(x)$ where $\sigma_i|F_0$ are different.  A little Galois theory shows that if $\theta$ is described as above and $\theta_1 \neq \pm \theta$ is an algebraic conjugate then $\theta_1^2=\tau(\ga)$ and $E=\Q[\theta+\theta_1]$.  One can see that the element $(\theta+\theta_1)^2$ is a primitive element of the real quadratic extension $\Q(\sqrt{N_{F_0|\Q}(\ga)})\neq F_0$ since $Gal(N|\Q)=D_{2,4}$.  The Lemma is thus established. $\blacksquare$

\end{caja}
We state the following proposition and prove only the cases not included in Moonen and Zarhin ~\cite{MZ}:
\begin{prop}\label{abe} Let $A_i$ be abelian varieties of dimension $1$ or $2$.  Then the Hodge conjecture holds for $A_1 \times \cdots \times A_r$ for $r$ an arbitrary natural number.
\end{prop}
\begin{caja}{Proof of Proposition \ref{abe}:} By Goursat's Lemma ~\cite{Lewis} ~\cite{Ribet} and the results of Hazama ~\cite{Hazama}  and Moonen-Zarhin ~\cite{MZ} one needs only prove the following statement.  For $A_i$ such that $\mbox{dim }A_i\leq 2$ and $\mbox{Hom}(A_1,A_2)=0$, one has $Hg(A_1 \times A_2)=Hg(A_1)\times Hg(A_2).$  By Hazama ~\cite{Lewis}, Moonen and Zarhin ~\cite{MZ} the only case left is the following.

  Let $A_i$ be simple abelian varieties of CM type.  Then $Hg(A_i)=U_{F_i}(1)$ and $Hg(A_1 \times A_2) \subset Hg(A_1)\times Hg(A_2)$ surjects onto both components, so either $Hg(A_1\times A_2)$ is simple (and the projections are isogenies) or the former inclusion is an equality.  Suppose that the projections are isogenies; in this case, $T(A_1)\otimes T(A_2)$ has a Hodge class (in fact, four such classes), and thus there is a Hodge isomorphism $T(A_1)\cong T(A_2)$.  This implies that $E_i=End_{hg}T(A_i)$ are isomorphic number fields; in the case where the Galois group of $N_1=F_1^{gal}$ over $\Q$ is $D_{2,4}$, we have $E_1\simeq E_2$ and it follows from Lemma \ref{diedral} that $F_1\simeq F_2$ as well.  The Proposition follows in this case from ~\cite{MZ} Proposition 4.2. For the remaining cases, there is only one CM type for $F$ up to automorphisms of $F|\Q$ and the proof is similar. $\blacksquare$

\end{caja}
  Now Theorem \ref{principal} follows easily from Proposition \ref{abe} and Remark \ref{ene}.

\section{The case of powers of a K3 surface}

Let $X$ be a K3 surface, and let ${\cal H}^\bullet(X) \subset H^\bullet(X)$ be the ring of Hodge classes of $X$.  Then $H^\bullet(X)=T(X)\oplus {\cal H}^\bullet(X).$  $T(X)$ is an irreducible Hodge structure ~\cite{vG} ~\cite{Zarhin}, and if $E=End_{hg}T(X)$ we have an inclusion $$E \hookrightarrow End_\C(H^{2,0}(X))=\C$$ which renders $E$ a number field.  It can be shown that $E$ is either totally real or CM ~\cite{Zarhin}.

The following proposition holds:

\begin{prop}\label{potencias} The Hodge conjecture for $X^n$, for arbitrary $n$, holds if it holds for $X \times X$.
\end{prop}
{\bf Proof:}  The ring of Hodge classes  ${\cal H}^\bullet(X^n)$ is, by the above, generated by the Hodge classes in the tensor powers of $T(X)$ up to order $n$ and by pullbacks of algebraic classes on $X$ via the canonical projections.  Thus our result amounts to show that the ring of tensor invariants of the $Hg(X)$-module $T(X)$ is generated by those of degree $2$ as an algebra; it is known (see ~\cite{Zarhin}) that $Hg(X)_\C$ is isomorphic to a product of special orthogonal or general linear groups, which shows (see ~\cite{Ribet}) that the ring of tensor invariants of $Hg(X)$ is generated by the degree-$2$ invariants, thereby establishing the result. $\blacksquare$

We now prove the Hodge conjecture for self-products of a K3 surface $X$ in the case where $E$ is a CM field.  We need the following elementary lemma.
\begin{lemma}Let $E$ be a CM number field.  Then $E$ is spanned as a vector space over $\Q$ by elements $\ga_i \in E$ such that $\ga_i \overline{\ga}_i=1$. \label{moduno}
\end{lemma}
{\bf Proof of Lemma \ref{moduno}:}  Let $\chi_0:E^\times \rightarrow E^\times$ be given by $\chi_0(\ga)=\ga/\overline{\ga}.$  Suppose that the images of $\chi_0$ do not span $E$ over $\Q$; then there exists $\theta \in E$ such that
$$\mbox{Tr}_{E|\Q}(\theta \chi_0(\ga))=0 \mbox{ for all } \ga \in E^\times.$$
Now let $\chi_\sigma= \sigma \circ \chi_0$ for $\sigma:E \hookrightarrow \C$ an embedding; by Artin's linear independence of characters, there are embeddings $\sigma \neq \tau$ such that $\chi_\sigma=\chi_\tau$, which amounts to saying that $\sigma(\ga)/\tau(\ga)$ is always real.  It is not difficult to see that, since $E$ is non-real CM, this cannot hold if $\sigma$ and $\tau$ are different; indeed, evaluating at $\ga$ and $1+\ga$ for $\ga \in E$ neither real nor purely imaginary, we see that $1+\sigma(\ga)$ does not belong to $\R\left(1+\tau(\ga)\right)$, which leads to a contradiction, thereby establishing the Lemma.
$\blacksquare$

  We are now ready to prove our Theorem.  See  Morrison ~\cite{M} for an earlier result in this direction.
\begin{theorem}\label{k3} Let $X$ be a K3 surface such that $E=End_{hg}T(X)$ is a CM field.  Then the Hodge conjecture holds for arbitrary powers of $X$.
\end{theorem}

\begin{caja}{\bf Proof of Theorem \ref{k3}:} By the above Lemma \ref{moduno}, it suffices to prove algebraicity for $\ga \in E$ such that $\ga\cdotp \overline{\ga}=1$, \textit{i.e.} for the Hodge isometries of the polarised Hodge structure $(T(X), Q)$ ~\cite{Zarhin}.  This is a result established by Mukai, by refining former results on his theory of moduli:\end{caja}
\begin{theorem}~\cite{Mukai}  Let $X_1$ and $X_2$ be K3 surfaces, and let $\psi:T(X_1) \rightarrow T(X_2)$ be a Hodge isometry.  Then $\psi$ is induced by an algebraic cycle.
\end{theorem}


\begin{thebibliography}{xxx}

\bibitem{Beauville} \rm{A. Beauville}, \textit{Complex Algebraic Surfaces}, 2nd English edition, Cambridge Univ. Press, 1996.

\bibitem{DeligneK3} \rm{P. Deligne}, \textit{La Conjecture de Weil pour les surfaces K3}, Inv. Math 15 (1972), 206-226.
\bibitem{DMOS} \rm{P. Deligne, J.S. Milne, A. Ogus, K.Y. Shih}, \textit{Hodge cycles, motives and Shimura varieties}, LNM 900.  Springer-Verlag, Berlin-Heidelberg-New York, 1982.
\bibitem{Freitag} \rm{E. Freitag}, \textit{\"Uber die Struktur der Funktionenk\"orper zu hyperabelschen Gruppen.I}. J.Reine und Angew. Math 247(1971), 97-117.
\bibitem{vG} \rm{B.L. van Geemen}, \textit{Kuga-Satake varieties and the Hodge conjecture}.  \textit{The Geometry and Arithmetic of Algebraic Cycles}. J.D. Lewis, B.Brent Gordon, S.M\"uller-Stach, N.Yui (eds.), NATO Science Series 548, World Scientific 2001.

\bibitem{Hazama} F. Hazama, \textit{Algebraic cycles on nonsimple abelian varieties}, Duke Math. J. Vol. 58 (1989) No. 1, pp.31-37.  
\bibitem{Jannsen} \rm{U. Jannsen}, \textit{Motives, numerical equivalence, and semisimplicity.}  Inv. Math. 107 (1992) 447-452.
\bibitem{Kleiman} S. Kleiman, \textit{Algebraic cycles and the Weil conjectures.} In: Dix Expos\'es sur la Cohomologie des Sch\'emas, North Holland, Amsterdam 1968, pp. 359-386.     
\bibitem{Lewis} \rm{J. D. Lewis,} \textit{A survey of the Hodge conjecture, }   2nd Edition, AMS Providence, RI 1998 (with an Appendix by B. Brent Gordon on the Hodge conjecture on abelian varieties.)
\bibitem{Manin} Yu I. Manin, \textit{Correspondences, motifs and monoidal transformations}, Math. USSR Sbornik 6 (1968), 439-470.
\bibitem{MZ} \rm{B.J.J. Moonen, Yu G. Zarhin}, \textit{Hodge cycles on abelian varieties of low dimension}.  Math.Ann. 315(1999), No. 4, 711-733.
\bibitem{M} D. R. Morrison, \textit{Algebraic cycles on products of surfaces.}  In: Proceedings, Algebraic Geometry Symposium, T\^ohoku Univ. 1984, pp. 194-210.  
\bibitem{Mukai1} \rm{S.Mukai}, \textit{On the moduli space of bundles of K3 surfaces, I.}  Vector Bundles on Algebraic Varieties, Bombay Colloquium 1984, Oxford 1987.
\bibitem{Mukai} \rm{S. Mukai}, \textit{Vector Bundles on a K3 Surface}. Proceedings of the ICM, Beijing 2002, Vol. II (invited lectures), pp. 485-502.  Higher Education Press -- World Scientific, Beijing 2002. 
\bibitem{MAV} \rm{D. Mumford}, Abelian Varieties, Oxford University Press, Oxford 1970.
\bibitem{Mumford} \rm{D. Mumford}, \textit{On Shimura's paper ``Discontinuous groups and abelian varieties"}.  Math. Ann. 181 (1969), 345-351.
\bibitem{Murre} \rm{J.P. Murre}, \textit{On the motive of an algebraic surface}, J. Reine und Angew. Math. 409(1990), 190-204.

\bibitem{T} \rm{J.J. Ram\'on Mar\'i,} \textit{On the Hodge conjecture for products of certain surfaces.}  PhD Thesis, U. Durham, UK 2003.
\bibitem{Ribet} \rm{K. Ribet}, \textit{Hodge classes on certain Abelian varieties}, Amer. J. Math 115 (1983), 523-538.
\bibitem{Schoen} \rm{C. Schoen}, \textit{Hodge classes on self-products of a variety with an automorphism}, Comp. Math. 65(1988) 3-32.
\bibitem{Scholl} \rm{A.J. Scholl}, \textit{Classical Motives}, Proc. Symp. Pure Math. 55 (Seattle, WA 1991), Part 1, pp. 163-187.  U. Jannsen, S. Kleiman, J. P. Serre (eds.)  AMS, Providence, RI 1995.
\bibitem{ST} G. Shimura, Y. Taniyama, \textit{Complex Multiplication of abelian varieties,} Publ. Math. Soc. Japan, vol. 6. Math. Soc. Japan, 1961. 
\bibitem{Serrano}  \rm{F. Serrano}, \textit{The Picard group of a quasi-bundle}, Manuscripta Mathematica 73 (1991), 63-82.
\bibitem{Ueno} \rm{K. Ueno}, \textit{Classification Theory of Algebraic Varieties and Compact Complex Spaces}, LNM 439, Springer-Verlag, Berlin-Heidelberg, 1975.
\bibitem{Zarhin} \rm{Yu. G. Zarhin},\textit{Hodge groups of K3 surfaces,} J. Reine und Angew. Math. 341 (1983), 193-220.
\end{thebibliography}
\end{document}